# Further bounds on $q$-numerical radius of Hilbert space operators

Satyajit Sahoo[a], Nirmal Chandra Rout[b]

[a]*Department of Mathematics, School of Basic Sciences, Indian Institute of Technology Bhubaneswar, Odisha 752050, India*
[b]*Department of Mathematics, Bhadrak Autonomous College, Bhadrak, Odisha 756100, India*

**Abstract**

In this article, we developed a series of new inequalities involving the $q$-numerical radius for operators and $2 \times 2$ operator matrices. These inequalities serve to establish both lower and upper bounds for the $q$-numerical radius of operators. Additionally, we established $q$-numerical radius inequalities for operators via Buzano inequality.

## 1. Introduction

The $q$-numerical radius is a generalization of the standard numerical radius of an operator. It was introduced in operator theory to study certain types of operator norms and spectral properties. The concept of the $q$-numerical radius is particularly used in the context of $q$-norms and more generalized settings of bounded linear operators on Hilbert spaces or Banach spaces. The value of $q$ influences how the numerical radius behaves. For $q > 1$, the value $w_q(T)$ tends to be larger compared to the standard numerical radius, while for $0 < q < 1$, the $q$-numerical radius tends to be smaller. The $q$-numerical radius provides a flexible way to study operator behavior under different norms or scaling of the inner product. It can be used to explore different spectral properties of operators, particularly when analyzing the operator in non-standard settings, such as for $q$-norms or other spaces that are not strictly Hilbert spaces.

---

*Nirmal Chandra Rout
  *Email addresses:* ssahoomath@gmail.com (Satyajit Sahoo), nrout89@gmail.com (Nirmal Chandra Rout)

Let $\mathcal{H}$ be a complex Hilbert space with inner product $\langle \cdot, \cdot \rangle$ and the corresponding norm $\|\cdot\|$. Let $\mathcal{L}(\mathcal{H})$ be the $C^*$-algebra of all bounded linear operators from $\mathcal{H}$ into itself. An operator $S \in \mathcal{L}(\mathcal{H})$ is said to be positive, and denoted $S \geq 0$, if $\langle Sx, x \rangle \geq 0$ for all $x \in \mathcal{H}$, and is called positive definite, denoted $S > 0$, if $\langle Sx, x \rangle > 0$ for all non zero vectors $x \in \mathcal{H}$. The *numerical range* of $S \in \mathcal{L}(\mathcal{H})$ is defined as $W(S) = \{\langle Sx, x \rangle : x \in \mathcal{H}, \|x\| = 1\}$ and the *numerical radius* of $S$, denoted by $w(S)$, is defined by $w(S) = \sup\{|z| : z \in W(S)\}$. It is known that the set $W(S)$ is a convex subset of the complex plane and that the numerical radius $w(\cdot)$ is a norm on $\mathcal{L}(\mathcal{H})$; being equivalent to the usual operator norm $\|S\| = \sup\{\|Sx\| : x \in \mathcal{H}, \|x\| = 1\}$. In fact, for every $S \in \mathcal{L}(\mathcal{H})$,

$$\frac{1}{2}\|S\| \leq w(S) \leq \|S\|. \tag{1.1}$$

The inequalities in (1.1) are sharp. If $S^2 = 0$, then the first inequality becomes an equality, on the other hand, the second inequality becomes an equality if $S$ is normal. In fact, for a nilpotent operator $T$ with $S^n = 0$, Haagerup and Harpe [20] showed that $w(S) \leq \|S\| \cos(\pi/(n+1))$. In particular, when $n = 2$, we get the reverse inequality of the first inequality in (1.1). The numerical radius has some significant properties, such as the power inequality:

$$w(S^n) \leq w^n(S) \quad \text{for} \quad n = 1, 2, \ldots. \tag{1.2}$$

For basic information about numerical radius, one can refer [19]. The author of [23, 25] improved the inequality (1.1) which is stated next. If $S \in \mathcal{L}(\mathcal{H})$, then

$$w(S) \leq \frac{1}{2}\||S| + |S^*|\| \leq \frac{1}{2}(\|S\| + \|S^2\|^{1/2}), \tag{1.3}$$

where $|S| = (S^*S)^{1/2}$ is the absolute value of $S$, and

$$\frac{1}{4}\|S^*S + SS^*\| \leq w^2(S) \leq \frac{1}{2}\|S^*S + SS^*\|. \tag{1.4}$$

The inequalities in (1.3) refines the second inequality in (1.1). For applications of these inequalities, one can refer [23, 24].

Dragomir [17] showed the following numerical radius inequality involving the product of two operators:

$$w^r(S^*T) \leq \frac{1}{2}\||T|^{2r} + |S|^{2r}\|, \quad r \geq 1. \tag{1.5}$$

The Schwarz inequality for positive operators reads that if $S$ is a positive operator in $\mathcal{L}(\mathcal{H})$, then

$$|\langle Sx, y \rangle|^2 \leq \langle Sx, x \rangle \langle Sy, y \rangle, \tag{1.6}$$

for any vectors $x, y \in \mathcal{H}$. In 1952, Kato [22] introduced a companion inequality of (1.6), called the mixed Schwarz inequality, which asserts

$$|\langle Sx, y\rangle|^2 \leq \langle |S|^{2\alpha} x, x\rangle \langle |S^*|^{2(1-\alpha)} y, y\rangle, \qquad 0 \leq \alpha \leq 1, \tag{1.7}$$

for all operators $S \in \mathcal{L}(\mathcal{H})$ and any vectors $x, y \in \mathcal{H}$, where $|S| = (S^*S)^{1/2}$. In particular, the following inequality

$$|\langle Sx, y\rangle| \leq \sqrt{\langle |S| x, x\rangle \langle |S^*| y, y\rangle} \quad (\text{see}[21, \text{pp } 75-76]). \tag{1.8}$$

If $\mathcal{H} = \mathcal{H}_1 \oplus \mathcal{H}_2$ and $T \in \mathcal{L}(\mathcal{H})$, then $T$ can be written as a block-matrix [27]

$$T = \begin{bmatrix} I_1^* T I_1 & I_1^* T I_2 \\ I_2^* T I_1 & I_2^* T I_2 \end{bmatrix}, \tag{1.9}$$

where $I_j \in \mathcal{L}(\mathcal{H}_j, \mathcal{H})$, and $I_j(x) = x$. If $T, S \in \mathcal{L}(\mathcal{H})$, then

$$w\left(\begin{bmatrix} T & O \\ O & S \end{bmatrix}\right) = \max\{w(T), w(S)\}, \tag{1.10}$$

and

$$\left\|\begin{bmatrix} T & O \\ O & S \end{bmatrix}\right\| = \left\|\begin{bmatrix} O & T \\ S & O \end{bmatrix}\right\| = \max\{\|T\|, \|S\|\}, \tag{1.11}$$

In a similar way, the $q$-numerical range is defined by

$$W_q(T) = \{\langle Tx, y\rangle : x, y \in \mathcal{H}, \|x\| = \|y\| = 1, \langle x, y\rangle = q\},$$

while the *q-numerical radius* is defined by

$$w_q(T) = \sup\{|z| : z \in W_q(T)\}. \tag{1.12}$$

One can observe that the $q$-numerical radius is a generalization of the classical numerical radius, for $|q| = 1$. This observation follows from the fact that equality would have to hold in the Cauchy-Schwarz inequality $|q| = |\langle x, y\rangle| \leq \|x\|\|y\| = 1$, provided $|q| = 1$. It is clear that $y = \lambda x$ would need to be true for some $\lambda \in \mathbb{C}, |\lambda| = 1$, hence $|\langle Sx, y\rangle| = |\langle Sx, x\rangle|$.

In 1977, Marcus and Andresen [5] introduced the notion of $q$-numerical range on $n$-dimensional unitary space with an inner product. They proved that the set generated by rotating $q$-numerical range about the origin is an annulus for $q \in \mathbb{C}$ with $q \leq 1$. They also established the inner and outer radii of this annulus for hermitian operators. In 1984, Nam-Kiu Tsing [6], established the convexity of the $q$-numerical range. In 1994, C. K. Li et al. [7] described some elementary properties for the $q$-numerical range on finite-dimensional spaces, see also [8]. In 2002, M.T. Chien and H. Nakazato [9] described the boundary of the $q$-numerical range of a square matrix using its Davis-Wielandt shell. In 2005, R. Rajic [10] considered a generalization of the $q$-numerical range. In 2007, M.T. Chien and H. Nakazato [11] again studied the $q$-numerical radius of weighted unilateral and bilateral shift operators and computed the $q$-numerical radius of shift operators with periodic weights. In 2012, M. T. Chien, [12] investigated the $q$-numerical radius of a weighted shift operator when its weights are in geometric sequence and periodic sequence. Recently, Moghaddam et al. [2, 3, 4] established some upper and lower bounds for the $q$-numerical radius of bounded linear operators which generalize some classical numerical radius inequalities. Also, they have presented some useful examples to compare the sharpness of their inequalities for different values of $q \in (0, 1)$. Very recently, Stankovic et al. [1] investigated certain properties for the $q$-numerical radius with a theoretical approach and presented an improved version of some earlier results from [2]. The authors of [1] improved the result [2, Theorem 2.1]. It states that for $T \in \mathcal{L}(\mathcal{H})$ and $q \in \overline{\mathbb{D}}$,

$$\frac{|q|}{2}\|T\| \leq w_q(T) \leq \|T\|. \tag{1.13}$$

Earlier results are proved for $q \in [0, 1]$, while in [1] some of the results are established with $q \in \overline{\mathbb{D}}$, the closed unit disc in $\mathbb{C}$. Also, they established several $q$-numerical radius inequalities for operator matrices defined on the direct sum of Hilbert spaces. Motivated by the work of several authors we have presented the $q$-numerical radius of some special type of operator matrices.

The objective of the article is to study the generalized numerical radius of operators and operator matrices and investigate certain analytical properties that may influence interested readers for further investigation. In this aspect, the rest of the article is organized as follows: In Section 2, we establish certain inequalities for the $q$-numerical radius for operators and operator matrices. Section 3 deals with several inequalities involving $q$-numerical radius for $2 \times 2$ operator matrices via Buzano inequality followed by comparisons of these results with the usual numerical radius, $\mathbb{A}$-numerical radius and Hilbert-Schmidt numerical radius. In order to make progress in our research work, we need the following lemmas to prove our results.

The very first lemma represents a consequence of the results obtained in [18, Proposition 3.1], as well as (1.12).

**Lemma 1.1.** *Let $T, S \in \mathcal{L}(\mathcal{H}), q \in \overline{\mathbb{D}}$ and $\lambda \in \mathbb{C}$. Then we have following properties:*

(i) $w_q(\lambda T) = |\lambda| w_q(T)$.

(ii) $w_q(T+S) \leq w_q(T) + w_q(S)$.

(iii) $w_q(U^*TU) = w_q(T)$, where $U \in \mathcal{B}(\mathcal{H})$ is an unitary operator.

(iv) $w_{\lambda q}(T) = w_q(T)$ for all $\lambda \in \mathbb{C}$ with $|\lambda| = 1$.

It is well known that, for $T_1, T_2 \in \mathcal{L}(H)$, $w\left(\begin{bmatrix} T_1 & O \\ O & T_2 \end{bmatrix}\right) = \max\{w(T_1), w(T_2)\}$, but the same result does not hold for $q$-numerical radius of the operator matrix [2]. In the forthcoming lemmas, we have collected some results for $q$-numerical radius of $2 \times 2$ and $n \times n$ operator matrices, which are essential to prove our main results.

**Lemma 1.2.** [1, Lemma 2.1] *Let $T_1, T_2 \in \mathcal{L}(\mathcal{H}), q \in \bar{\mathbb{D}}$ and $\theta \in \mathbb{R}$. Then*

(i) $w_q\left(\begin{bmatrix} O & T_1 \\ T_2 & O \end{bmatrix}\right) = w_q\left(\begin{bmatrix} O & T_2 \\ T_1 & O \end{bmatrix}\right)$.

(ii) $w_q\left(\begin{bmatrix} O & T_1 \\ e^{i\theta}T_2 & O \end{bmatrix}\right) = w_q\left(\begin{bmatrix} O & T_1 \\ T_2 & O \end{bmatrix}\right)$ *for any $\theta \in \mathbb{R}$.*

(iii) $w_q\left(\begin{bmatrix} T_1 & O \\ O & T_2 \end{bmatrix}\right) = w_q\left(\begin{bmatrix} T_2 & O \\ O & T_1 \end{bmatrix}\right)$.

The following Lemma follows directly from [28, p 107].

**Lemma 1.3.** *Let $T_1, T_2, T_3, T_4 \in \mathcal{L}(\mathcal{H})$. Then*

(i) $w_q\left(\begin{bmatrix} T_1 & O \\ O & T_4 \end{bmatrix}\right) \leq w_q\left(\begin{bmatrix} T_1 & T_2 \\ T_3 & T_4 \end{bmatrix}\right)$.

(ii) $w_q\left(\begin{bmatrix} O & T_2 \\ T_3 & O \end{bmatrix}\right) \leq w_q\left(\begin{bmatrix} T_1 & T_2 \\ T_3 & T_4 \end{bmatrix}\right)$.

The following result is from [1, Theorem 1.5].

**Lemma 1.4.** [1, Theorem 1.5] *Let $(\mathcal{H}_n)_{n \in \mathbb{N}}$ be a sequence of Hilbert spaces and let $T_n \in \mathcal{L}(\mathcal{H}_n)$ for all $n \in \mathbb{N}$. If $q \in \bar{\mathbb{D}} \smallsetminus \{0\}$. Then*

$$\sup_{n \in \mathbb{N}} w_q(T_n) \leq w_q\left(\bigoplus_{n=1}^{+\infty} T_n\right) \leq \frac{|q| + 2\sqrt{1 - |q|^2}}{|q|} \sup_{n \in \mathbb{N}} w_q(T_n).$$

As a special case of the above result, we have the following.

**Lemma 1.5.** [1, Corollary 5.1] *Let $T, S \in \mathcal{L}(\mathcal{H})$ and $q \in (0, 1]$. Then*

$$\max\{w_q(T), w_q(S)\} \leq w_q\left(\begin{bmatrix} T & 0 \\ 0 & S \end{bmatrix}\right) \leq \frac{q + 2\sqrt{1 - q^2}}{q} \max\{w_q(T), w_q(S)\}.$$

**Lemma 1.6.** *Let $T \in \mathcal{L}(\mathcal{H})$ and $q \in \overline{\mathbb{D}}$. Then*

$$|q|^{n-1} w_q(T^n) \leq w_q^n(T) \quad \text{for all } n \in \mathbb{N}. \tag{1.14}$$

**Lemma 1.7.** [26] *Let $0 \leq q \leq 1$ and $T \in \mathcal{M}_2(\mathbb{C})$. Then $T$ is unitarily similar to $e^{it}\begin{bmatrix} \gamma & a \\ b & \gamma \end{bmatrix}$ for some $0 \leq t \leq 2\pi$ and $0 \leq b \leq a$. Also,*

$$W_q(T) = e^{it}\{\gamma q + r((c+pd)\cos(s) + i(d+pc)\sin(s)) : 0 \leq r \leq 1, 0 \leq s \leq 2\pi\},$$

*with $c = \frac{a+b}{2}$, $d = \frac{a-b}{2}$ and $p = \sqrt{1-q^2}$.*

## 2. $q$-numerical radius inequalities for operators and operator matrices

In this section, we have established several results on $q$-numerical radius of operators and $2 \times 2$ operator matrices. Some special cases of our results lead to a comparison with some of the results in the literature. For better understanding, we have interpreted the geometry of the $q$-numerical range with a few interesting examples for interested readers. The first result of this section provides an upper bound for $q$-numerical radius of an operator.

**Theorem 2.1.** *Let $S \in \mathcal{L}(\mathcal{H})$ and $0 < q < 1, 0 \leq \alpha \leq 1$ and $r \geq 1$, then*

$$w_q^{2r}(S) \leq q^2 \left\| \alpha |S|^{2r} + (1-\alpha)|S^*|^{2r} \right\| + (1-q^2)(\alpha \| |S|^{2r} \| + (1-\alpha) \| |S^*|^{2r} \|)$$
$$+ 2q\sqrt{1-q^2} \left\| |S^*|^{2r(1-\alpha)} \right\| \left\| |S|^{2r\alpha} \right\|.$$

*Proof.*

$$|\langle Sx, y\rangle|^{2r} \leq \left(\langle |S|^{2\alpha}x, x\rangle \langle |S^*|^{2(1-\alpha)}y, y\rangle\right)^r$$
$$= \langle |S|^{2\alpha}x, x\rangle^r \langle |S^*|^{2(1-\alpha)}y, y\rangle^r$$
$$= \langle |S|^{2\alpha r}x, x\rangle \langle |S^*|^{2(1-\alpha)r}y, y\rangle$$
$$= \langle |S|^{2\alpha r}x, x\rangle \langle |S^*|^{2(1-\alpha)r}(qx + \sqrt{1-q^2}z), qx + \sqrt{1-q^2}z\rangle$$
$$= \langle |S|^{2\alpha r}x, x\rangle \Big[\langle |S^*|^{2(1-\alpha)r}qx, qx\rangle + \langle |S^*|^{2(1-\alpha)r}qx, \sqrt{1-q^2}z\rangle$$
$$+ \langle |S^*|^{2(1-\alpha)r}\sqrt{1-q^2}z, qx\rangle + \langle |S^*|^{2(1-\alpha)r}\sqrt{1-q^2}z, \sqrt{1-q^2}z\rangle\Big]$$

$$= \langle |S|^{2\alpha r} x, x \rangle \Big[ q^2 \langle |S^*|^{2(1-\alpha)r} x, x \rangle + q\sqrt{1-q^2} \langle |S^*|^{2(1-\alpha)r} x, z \rangle$$
$$+ q\sqrt{1-q^2} \langle |S^*|^{2(1-\alpha)r} z, x \rangle + (1-q^2) \langle |S^*|^{2(1-\alpha)r} z, z \rangle \Big]$$
$$= q^2 \langle |S|^{2\alpha r} x, x \rangle \langle |S^*|^{2(1-\alpha)r} x, x \rangle + q\sqrt{1-q^2} \Big( \langle |S^*|^{2(1-\alpha)r} x, z \rangle$$
$$+ \langle |S^*|^{2(1-\alpha)r} z, x \rangle \Big) \langle |S|^{2\alpha r} x, x \rangle + (1-q^2) \langle |S|^{2\alpha r} x, x \rangle \langle |S^*|^{2(1-\alpha)r} z, z \rangle$$
$$\leq q^2 \langle |S|^{2r} x, x \rangle^{\alpha} \langle |S^*|^{2r} x, x \rangle^{(1-\alpha)} + (1-q^2) \langle |S|^{2r} x, x \rangle^{\alpha} \langle |S^*|^{2r} z, z \rangle^{(1-\alpha)}$$
$$+ q\sqrt{1-q^2} 2\Re \Big( \langle |S^*|^{2(1-\alpha)r} x, z \rangle \Big) \langle |S|^{2\alpha r} x, x \rangle$$
$$\leq q^2 \Big( \alpha \langle |S|^{2r} x, x \rangle + (1-\alpha) \langle |S^*|^{2r} x, x \rangle \Big) + (1-q^2) \Big( \alpha \langle |S|^{2r} x, x \rangle + (1-\alpha) \langle |S^*|^{2r} z, z \rangle \Big)$$
$$+ 2q\sqrt{1-q^2} \Big| \langle |S^*|^{2(1-\alpha)r} x, z \rangle \Big| \langle |S|^{2\alpha r} x, x \rangle$$
$$\leq q^2 \langle (\alpha |S|^{2r} + (1-\alpha) |S^*|^{2r}) x, x \rangle + (1-q^2) \Big( \alpha \langle |S|^{2r} x, x \rangle + (1-\alpha) \langle |S^*|^{2r} z, z \rangle \Big)$$
$$+ 2q\sqrt{1-q^2} \, \||S^*|^{2(1-\alpha)r}\| \, \||S|^{2\alpha r}\| \, \|x\|^3 \|z\|$$
$$\leq q^2 \, \|\alpha |S|^{2r} + (1-\alpha) |S^*|^{2r}\| + (1-q^2)(\alpha \||S|^{2r}\| + (1-\alpha) \||S^*|^{2r}\|)$$
$$+ 2q\sqrt{1-q^2} \, \||S^*|^{2(1-\alpha)r}\| \, \||S|^{2\alpha r}\|.$$

So, taking supremum over all $x \in \mathcal{H}$ with $\|x\| = 1 = \|y\|$ and $\langle x, y \rangle = q$, we find that

$$w_q^{2r}(S) \leq q^2 \, \|\alpha|S|^{2r} + (1-\alpha)|S^*|^{2r}\| + (1-q^2)(\alpha\||S|^{2r}\| + (1-\alpha)\||S^*|^{2r}\|) + 2q\sqrt{1-q^2} \, \||S^*|^{2(1-\alpha)r}\| \, \||S|^{2\alpha r}\|.$$

□

The following remark can be obtained by taking $r = 1$.

**Remark 2.2.** Let $S \in \mathcal{L}(\mathcal{H})$ and $0 < q < 1, 0 \leq \alpha \leq 1$, then

$$w_q^2(S) \leq q^2 \, \|\alpha|S|^2 + (1-\alpha)|S^*|^2\| + (1-q^2)(\alpha\||S|^2\| + (1-\alpha)\||S^*|^2\|) + 2q\sqrt{1-q^2} \, \||S^*|^{2(1-\alpha)}\| \, \||S|^{2\alpha}\|. \tag{2.1}$$

**Remark 2.3.** For $q \to 1$, we have $w^2(S) \leq \|\alpha|S|^2 + (1-\alpha)|S^*|^2\|$, which is a generalized version of right-hand side of the Kittaneh inequality (1.4). For $q \to 1$ and $\alpha = \frac{1}{2}$, we have $w^2(S) \leq \frac{1}{2} \||S|^2 + |S^*|^2\|$ which is the right-hand side of the inequality (1.4).

The following result is known from [1, Lemma 5.3].

**Lemma 2.4.** *Let $T, S \in \mathcal{L}(\mathcal{H})$, $q \in (0,1]$. Then*

$$\max\{w_q(T+S), w_q(T-S)\} \leq w_q\left(\begin{bmatrix} T & S \\ S & T \end{bmatrix}\right) \leq \frac{q + 2\sqrt{1-q^2}}{q} \max\{w_q(T+S), w_q(T-S)\}.$$

*In particular,*

$$w_q(S) \leq w_q\left(\begin{bmatrix} O & S \\ S & O \end{bmatrix}\right) \leq \frac{q + 2\sqrt{1-q^2}}{q} w_q(S).$$

We have a new upper bound for $q$-numerical radius.

**Lemma 2.5.** *Let $P, Q \in \mathcal{L}(\mathcal{H})$, $q \in (0,1]$. Then*

$$\max\left\{w_q(P+iQ), w_q(P-iQ)\right\} \leq w_q\left(\begin{bmatrix} Q & -P \\ P & Q \end{bmatrix}\right) \leq \frac{q + 2\sqrt{1-q^2}}{q} \max\left\{w_q(P+iQ), w_q(P-iQ)\right\}.$$

*Proof.* Let $T = \begin{bmatrix} iQ & -P \\ P & iQ \end{bmatrix}$ and $U = \frac{1}{\sqrt{2}}\begin{bmatrix} I & iI \\ iI & I \end{bmatrix}$. So, $U^* = \frac{1}{\sqrt{2}}\begin{bmatrix} I & -iI \\ -iI & I \end{bmatrix}$. It is not difficult to show that $U$ is a unitary operator on $\mathcal{H} \oplus \mathcal{H}$. Then $U^*TU = \begin{bmatrix} -i(P-Q) & O \\ O & i(P+Q) \end{bmatrix}$. Using the fact that $w_q(T) = w_q(U^*TU)$, we have

$$w_q(T) = w_q(U^*TU) = w_q\left(\begin{bmatrix} -i(P-Q) & O \\ O & i(P+Q) \end{bmatrix}\right).$$

Now, applying Lemma 1.5, we get

$$\max\{w_q(-i(P-Q)), w_q(i(P+Q))\} \leq w_q\left(\begin{bmatrix} iQ & -P \\ P & iQ \end{bmatrix}\right) \leq \frac{q + 2\sqrt{1-q^2}}{q} \max\{w_q(-i(P-Q)), w_q(i(P+Q))\}$$

$$= \max\{w_q(P-Q), w_q(P+Q)\} \leq w_q\left(\begin{bmatrix} iQ & -P \\ P & iQ \end{bmatrix}\right) \leq \frac{q + 2\sqrt{1-q^2}}{q} \max\{w_q(P-Q), w_q(P+Q)\}.$$

Replacing $Q$ by $-iQ$ in the above inequality, we have

$$\max\left\{w_q(P+iQ), w_q(P-iQ)\right\} \leq w_q\left(\begin{bmatrix} Q & -P \\ P & Q \end{bmatrix}\right) \leq \frac{q + 2\sqrt{1-q^2}}{q} \max\left\{w_q(P+iQ), w_q(P-iQ)\right\}.$$

□

**Remark 2.6.** *For $q = 1$, in Lemma 2.5, we have the equality for usual numerical radius see [16], and for $\mathbb{A}$-numerical radius and Hilbert-Schmidt numerical radius version of Lemma 2.5 one can see [13, 14].*

**Theorem 2.1.** *Let $P, Q, R, S \in \mathcal{L}(\mathcal{H})$, $q \in (0, 1]$. Then*

$$w_q\left(\begin{bmatrix} P & Q \\ R & S \end{bmatrix}\right) \leq \frac{q + 2\sqrt{1-q^2}}{2q}\Big\{\max\Big\{w_q((R-Q) + i(P+S)), w_q((R-Q) - i(P+S))\Big\}$$
$$+ \big(w_q(Q+R) + w_q(S-P)\big)\Big\}.$$

*Proof.* Let $U = \frac{1}{\sqrt{2}}\begin{bmatrix} I & -I \\ I & I \end{bmatrix}$. It can be shown that $U$ is a unitary operator on $\mathcal{H} \oplus \mathcal{H}$. Using the identity $w_q(T) = w_q(U^*TU)$, we have

$$w_q\left(\begin{bmatrix} P & Q \\ R & S \end{bmatrix}\right) = w_q\left(U^*\begin{bmatrix} P & Q \\ R & S \end{bmatrix}U\right)$$

$$= \frac{1}{2}w_q\left(\begin{bmatrix} P+Q+R+S & -P+Q-R+S \\ -P-Q+R+S & P-Q-R+S \end{bmatrix}\right)$$

$$= \frac{1}{2}w_q\left(\begin{bmatrix} P+S & Q-R \\ R-Q & P+S \end{bmatrix} + \begin{bmatrix} Q+R & S-P \\ S-P & -R-Q \end{bmatrix}\right)$$

$$\leq \frac{1}{2}\left\{w_q\left(\begin{bmatrix} P+S & -(R-Q) \\ R-Q & P+S \end{bmatrix}\right) + w_q\left(\begin{bmatrix} Q+R & S-P \\ S-P & -(Q+R) \end{bmatrix}\right)\right\}$$

$$\leq \frac{q + 2\sqrt{1-q^2}}{2q}\Big\{\max\Big\{w_q((R-Q) + i(P+S)), w_q((R-Q) - i(P+S))\Big\}$$
$$+ \big(w_q(Q+R) + w_q(S-P)\big)\Big\},$$

where we have used Lemma 2.5, Lemma 1.5 and Lemma 2.4 to obtain the last inequality. This completes the proof. $\square$

**Remark 2.7.** *If we take $P = S$, in the previous theorem, we have*

$$w_q\left(\begin{bmatrix} P & Q \\ R & P \end{bmatrix}\right) \leq \frac{q + 2\sqrt{1-q^2}}{2q}\left\{\max\{w_q((R-Q) + 2iP), w_q((R-Q) - 2iP)\} + w_q(Q+R)\right\}.$$

**Remark 2.8.** *For $q \to 1$, in Theorem 2.1 we have the usual numerical radius see [15], and for $\mathbb{A}$-numerical radius and Hilbert-Schmidt numerical radius version of Theorem 2.1 one can see [13, 14].*

In order to prove our result the following result from [2, Theorem 2.5] is essential for our purpose. If $T \in \mathcal{L}(\mathcal{H})$ with $T^2 = 0$, then for any $q \in [0, 1)$,

$$w_q(T) \le \sqrt{1 - \frac{3q^2}{4} + q\sqrt{1-q^2}}\|T\|. \qquad (2.2)$$

Since $\begin{bmatrix} T & T \\ -T & -T \end{bmatrix}^2 = \begin{bmatrix} O & O \\ O & O \end{bmatrix}$, so by (2.2)

$$w_q\left(\begin{bmatrix} T & T \\ -T & -T \end{bmatrix}\right) \le \sqrt{1 - \frac{3q^2}{4} + q\sqrt{1-q^2}} \left\| \begin{bmatrix} T & T \\ -T & -T \end{bmatrix} \right\| = 2\sqrt{1 - \frac{3q^2}{4} + q\sqrt{1-q^2}}\|T\|. \qquad (2.3)$$

**Theorem 2.9.** *Let $T_2, T_3 \in \mathcal{L}(\mathcal{H})$ and $q \in [0, 1)$. Then*

$$w_q\left(\begin{bmatrix} O & T_2 \\ T_3 & O \end{bmatrix}\right) \le \frac{q + 2\sqrt{1-q^2}}{q} \min\{w_q(T_2), w_q(T_3)\} + \left(\sqrt{1 - \frac{3q^2}{4} + q\sqrt{1-q^2}}\right) \min\{\|T_2 + T_3\|, \|T_2 - T_3\|\}.$$

*Proof.* Let $U = \frac{1}{\sqrt{2}} \begin{bmatrix} I & -I \\ I & I \end{bmatrix}$. It is clear that $U$ is a unitary operator. Using the identity $w_q(T) = w_q(U^*TU)$, we have

$$w_q\left(\begin{bmatrix} O & T_2 \\ T_3 & O \end{bmatrix}\right) = w_q\left(U^* \begin{bmatrix} O & T_2 \\ T_3 & O \end{bmatrix} U\right)$$

$$= \frac{1}{2} w_q\left(\begin{bmatrix} T_2 + T_3 & T_2 - T_3 \\ -(T_2 - T_3) & -(T_2 + T_3) \end{bmatrix}\right)$$

$$= \frac{1}{2} w_q\left(\begin{bmatrix} T_2 + T_3 & T_2 + T_3 \\ -(T_2 + T_3) & -(T_2 + T_3) \end{bmatrix} + \begin{bmatrix} O & -2T_3 \\ 2T_3 & O \end{bmatrix}\right)$$

$$\le \frac{1}{2} \left\{ w_q\left(\begin{bmatrix} T_2 + T_3 & T_2 + T_3 \\ -(T_2 + T_3) & -(T_2 + T_3) \end{bmatrix}\right) + w_q\left(\begin{bmatrix} O & -2T_3 \\ 2T_3 & O \end{bmatrix}\right) \right\}.$$

Now, using inequality (2.3) and Lemma 2.4, we have

$$w_q\left(\begin{bmatrix} O & T_2 \\ T_3 & O \end{bmatrix}\right) \le \left(\sqrt{1 - \frac{3q^2}{4} + q\sqrt{1-q^2}}\right) \|T_2 + T_3\| + \frac{q + 2\sqrt{1-q^2}}{q} w_q(T_3). \qquad (2.4)$$

Replacing $T_3$ by $-T_3$ in the inequality (2.4) and using Lemma 1.2, we get

$$w_q\left(\begin{bmatrix} O & T_2 \\ T_3 & O \end{bmatrix}\right) \le \left(\sqrt{1 - \frac{3q^2}{4} + q\sqrt{1-q^2}}\right) \|T_2 - T_3\| + \frac{q + 2\sqrt{1-q^2}}{q} w_q(T_3). \qquad (2.5)$$

From the inequalities (2.4) and (2.5), we have

$$w_q\left(\begin{bmatrix} O & T_2 \\ T_3 & O \end{bmatrix}\right) \leq \frac{q + 2\sqrt{1-q^2}}{q} w_q(T_3) + \left(\sqrt{1 - \frac{3q^2}{4}} + q\sqrt{1-q^2}\right) \min\{\|T_2 + T_3\|, \|T_2 - T_3\|\}. \tag{2.6}$$

Again, in the inequality (2.6), interchanging $T_2$ and $T_3$ and using Lemma 1.2(ii), we get

$$w_q\left(\begin{bmatrix} O & T_2 \\ T_3 & O \end{bmatrix}\right) \leq \frac{q + 2\sqrt{1-q^2}}{q} w_q(T_2) + \left(\sqrt{1 - \frac{3q^2}{4}} + q\sqrt{1-q^2}\right) \min\{\|T_2 + T_3\|, \|T_2 - T_3\|\}. \tag{2.7}$$

From the inequalities (2.6) and (2.7), we get

$$w_q\left(\begin{bmatrix} O & T_2 \\ T_3 & O \end{bmatrix}\right) \leq \frac{q + 2\sqrt{1-q^2}}{q} \min\{w_q(T_2), w_q(T_3)\} + \left(\sqrt{1 - \frac{3q^2}{4}} + q\sqrt{1-q^2}\right) \min\{\|T_2 + T_3\|, \|T_2 - T_3\|\}.$$

This completes the proof. □

**Remark 2.10.** *If we set $q \to 1$, then we get a result by Hirzallah et al. [15] and for $\mathbb{A}$-numerical radius and Hilbert-Schmidt numerical radius version of the above theorem one can see [13, 14].*

The next results deal with the lower bounds for $q$-numerical radius of certain $2\times 2$ operator matrices.

**Theorem 2.11.** *Let $T_2, T_3 \in \mathcal{L}(\mathcal{H})$, $q \in (0, 1]$. Then*

$$w_q\left(\begin{bmatrix} O & T_2 \\ T_3 & O \end{bmatrix}\right) \geq \max\{w_q(T_2), w_q(T_3)\} - \sqrt{1 - \frac{3q^2}{4}} + q\sqrt{1-q^2} \min\{\|T_2 + T_3\|, \|T_2 - T_3\|\}.$$

*Proof.* Let $U = \frac{1}{\sqrt{2}}\begin{bmatrix} I & -I \\ I & I \end{bmatrix}$. It can be shown that $U$ is $\mathbb{A}$-unitary. Then

$$\frac{1}{2}\begin{bmatrix} T_2 + T_3 & T_2 + T_3 \\ -(T_2 + T_3) & -(T_2 + T_3) \end{bmatrix} = U^* \begin{bmatrix} O & T_2 \\ T_3 & O \end{bmatrix} U - \begin{bmatrix} O & -T_3 \\ T_3 & O \end{bmatrix}. \tag{2.8}$$

So,

$$\begin{bmatrix} O & -T_3 \\ T_3 & O \end{bmatrix} = U^* \begin{bmatrix} O & T_2 \\ T_3 & O \end{bmatrix} U - \frac{1}{2}\begin{bmatrix} T_2 + T_3 & T_2 + T_3 \\ -(T_2 + T_3) & -(T_2 + T_3) \end{bmatrix}. \tag{2.9}$$

This implies

$$w_q\left(\begin{bmatrix} O & -T_3 \\ T_3 & O \end{bmatrix}\right) \leq w_q\left(U^*\begin{bmatrix} O & T_2 \\ T_3 & O \end{bmatrix}U\right) + \frac{1}{2}w_q\left(\begin{bmatrix} T_2+T_3 & T_2+T_3 \\ -(T_2+T_3) & -(T_2+T_3) \end{bmatrix}\right).$$

Which in turn implies that

$$w_q\left(\begin{bmatrix} O & -T_3 \\ T_3 & O \end{bmatrix}\right) \leq w_q\left(\begin{bmatrix} O & T_2 \\ T_3 & O \end{bmatrix}\right) + \frac{1}{2}w_q\left(\begin{bmatrix} T_2+T_3 & T_2+T_3 \\ -(T_2+T_3) & -(T_2+T_3) \end{bmatrix}\right).$$

Thus, using inequality (2.3) and Lemma 2.5

$$w_q(T_3) \leq w_q\left(\begin{bmatrix} O & T_2 \\ T_3 & O \end{bmatrix}\right) + \sqrt{1 - \frac{3q^2}{4}} + q\sqrt{1-q^2}\|T_2+T_3\|. \tag{2.10}$$

Replacing $T_3$ by $-T_3$ in the inequality (2.10) we have

$$w_q(T_3) \leq w_q\left(\begin{bmatrix} O & T_2 \\ T_3 & O \end{bmatrix}\right) + \sqrt{1 - \frac{3q^2}{4}} + q\sqrt{1-q^2}\|T_2-T_3\|. \tag{2.11}$$

Now, from inequality (2.10) and (2.11) that

$$w_q(T_3) \leq w_q\left(\begin{bmatrix} O & T_2 \\ T_3 & O \end{bmatrix}\right) + \sqrt{1 - \frac{3q^2}{4}} + q\sqrt{1-q^2}\min\left\{\|T_2+T_3\|, \|T_2-T_3\|\right\}. \tag{2.12}$$

Interchanging $T_2$ and $T_3$ in the inequality (2.12), we get

$$w_q(T_2) \leq w_q\left(\begin{bmatrix} O & T_2 \\ T_3 & O \end{bmatrix}\right) + \sqrt{1 - \frac{3q^2}{4}} + q\sqrt{1-q^2}\min\left\{\|T_2+T_3\|, \|T_2-T_3\|\right\}. \tag{2.13}$$

From inequalities (2.12) and (2.13), we have

$$\max\{w_q(T_2), w_q(T_3)\} \leq w_q\left(\begin{bmatrix} O & T_2 \\ T_3 & O \end{bmatrix}\right) + \sqrt{1 - \frac{3q^2}{4}} + q\sqrt{1-q^2}\min\left\{\|T_2+T_3\|, \|T_2-T_3\|\right\}. \tag{2.14}$$

Which proves the first inequality. $\square$

**Theorem 2.12.** *Let $T_2, T_3 \in \mathcal{L}(\mathcal{H})$. Then*

$$w_q^2\left(\begin{bmatrix} O & T_2 \\ T_3 & O \end{bmatrix}\right) \geq \frac{q}{2}\max\left\{w_q(T_2T_3+T_3T_2), w_q(T_2T_3-T_3T_2)\right\}.$$

*Proof.* Let us consider an unitary operator $U = \begin{bmatrix} O & I \\ I & O \end{bmatrix}; T = \begin{bmatrix} O & T_2 \\ T_3 & O \end{bmatrix}$. Now,

$$T^2 + (U^*TU)^2 = \begin{bmatrix} T_2T_3 + T_3T_2 & O \\ O & T_3T_2 + T_2T_3 \end{bmatrix}.$$

Hence by using Lemma 1.5 we obtain

$$w_q(T_2T_3 + T_3T_2) \leq w_q\left(\begin{bmatrix} T_2T_3 + T_3T_2 & O \\ O & T_3T_2 + T_2T_3 \end{bmatrix}\right)$$
$$= w_q\left(T^2 + (U^*TU)^2\right)$$
$$\leq w_q\left(T^2\right) + w_q\left((U^*TU)^2\right)$$
$$\leq \frac{1}{q}w_q^2(T) + \frac{1}{q}w_q^2(U^*TU)$$
$$= \frac{1}{q}(w_q^2(T) + w_q^2(T))$$
$$= \frac{2}{q}w_q^2(T).$$

So,

$$w_q(T_2T_3 + T_3T_2) \leq \frac{2}{q}w_q^2(T). \tag{2.15}$$

Using similar argument to $(T)^2 - (U^*TU)^2$, we have

$$w_q(T_2T_3 - T_3T_2) \leq \frac{2}{q}w_q^2(T). \tag{2.16}$$

Combining (2.15) and (2.16) we get

$$w_q^2\left(\begin{bmatrix} O & T_2 \\ T_3 & O \end{bmatrix}\right) \geq \frac{q}{2}\max\left\{w_q(T_2T_3 + T_3T_2), w_q(T_2T_3 - T_3T_2)\right\}.$$

□

**Corollary 2.13.** *Let $T_1, T_2, T_3, T_4 \in \mathcal{L}(\mathcal{H})$. Then*

$$w_q\left(\begin{bmatrix} T_1 & T_2 \\ T_3 & T_4 \end{bmatrix}\right) \geq \max\left\{w_q(T_1), w_q(T_4), \sqrt{\frac{q}{2}}\left(w_q(T_2T_3 + T_3T_2)\right)^{\frac{1}{2}}, \sqrt{\frac{q}{2}}\left(w_q(T_2T_3 - T_3T_2)\right)^{\frac{1}{2}}\right\}.$$

*Proof.* Based on Lemma 1.3, Lemma 1.5 and Theorem 2.12 we have

$$w_q\left(\begin{bmatrix} T_1 & T_2 \\ T_3 & T_4 \end{bmatrix}\right) \geq \max\left\{w_q\left(\begin{bmatrix} T_1 & O \\ O & T_4 \end{bmatrix}\right), w_q\left(\begin{bmatrix} O & T_2 \\ T_3 & O \end{bmatrix}\right)\right\}$$

$$\geq \max\left\{w_q(T_1), w_q(T_4), \sqrt{\frac{q}{2}}\left(w_q(T_2T_3 + T_3T_2)\right)^{\frac{1}{2}}, \sqrt{\frac{q}{2}}\left(w_q(T_2T_3 - T_3T_2)\right)^{\frac{1}{2}}\right\}.$$

□

**Theorem 2.14.** *Let $T_2, T_3 \in \mathcal{L}(\mathcal{H})$. Then for $n \in \mathbb{N}$*

$$w_q\left(\begin{bmatrix} O & T_2 \\ T_3 & O \end{bmatrix}\right) \geq q^{2n-1}\left[\max\{w_q((T_2T_3)^n), w_q((T_3T_2)^n)\}\right]^{\frac{1}{2n}}. \tag{2.17}$$

*Proof.* Let $T = \begin{bmatrix} O & T_2 \\ T_3 & O \end{bmatrix}$. Then for $n \in \mathbb{N}$, $T^{2n} = \begin{bmatrix} (T_2T_3)^n & O \\ O & (T_3T_2)^n \end{bmatrix}$ and using Lemma 1.5 we obtain

$$\max\{w_q((T_2T_3)^n), w_q((T_3T_2)^n)\} \leq w_q\left(\begin{bmatrix} (T_2T_3)^n & O \\ O & (T_3T_2)^n \end{bmatrix}\right)$$

$$= w_q(T^{2n})$$

$$\leq \frac{1}{q^{2n-1}}w_q^{2n}(T) \quad \text{(by inequality (1.14))}$$

$$= \frac{1}{q^{2n-1}}w_q^{2n}\left(\begin{bmatrix} O & T_2 \\ T_3 & O \end{bmatrix}\right).$$

□

The following lemma is already proved by Hirzallah et al. [15] for the case of Hilbert space operators. Now we state here the result without proof for our purpose.

**Lemma 2.15.** *Let $T = \begin{bmatrix} T_1 & T_2 \\ T_2 & T_1 \end{bmatrix} \in \mathcal{L}(\mathcal{H} \oplus \mathcal{H})$ and $n \in \mathbb{N}$. Then $T^n = \begin{bmatrix} P & Q \\ Q & P \end{bmatrix}$ for some $P, Q \in \mathcal{L}(\mathcal{H})$ such that $P + Q = (T_1 + T_2)^n$ and $P - Q = (T_1 - T_2)^n$.*

The forthcoming result is analogous to Theorem 2.14

**Theorem 2.16.** *Let $T_1, T_2 \in \mathcal{L}(\mathcal{H})$. Then*

$$w_q\left(\begin{bmatrix} T_1 & T_2 \\ -T_2 & -T_1 \end{bmatrix}^{2n}\right) \leq \left[\max\{w_q(((T_1-T_2)(T_1+T_2))^n), w_q(((T_1+T_2)(T_1-T_2))^n)\}\right]$$

*for $n \in \mathbb{N}$.*

*Proof.* Let $T = \begin{bmatrix} T_1 & T_2 \\ -T_2 & -T_1 \end{bmatrix}$ and $R = T^2 = \begin{bmatrix} T_1^2 - T_2^2 & T_1T_2 - T_2T_1 \\ T_1T_2 - T_2T_1 & T_1^2 - T_2^2 \end{bmatrix}$. Using Lemma 2.15 we have there exist $P, Q \in \mathcal{L}(\mathcal{H})$ such that $R^n = \begin{bmatrix} P & Q \\ Q & P \end{bmatrix}$ with $P + Q = ((T_1^2 - T_2^2) + (T_1T_2 - T_2T_1))^n$ and $P - Q = ((T_1^2 - T_2^2) - (T_1T_2 - T_2T_1))^n$. So, $T^{2n} = \begin{bmatrix} P & Q \\ Q & P \end{bmatrix}$ with $P + Q = ((T_1 - T_2)(T_1 + T_2))^n$ and $P - Q = ((T_1 + T_2)(T_1 - T_2))^n$. Now, we have

$$w_q(T^{2n}) = w_q\left(\begin{bmatrix} P & Q \\ Q & P \end{bmatrix}\right)$$
$$\leq \frac{q + 2\sqrt{1-q^2}}{q} \max\{w_q(P+Q), w_q(P-Q)\} \quad \text{(by Lemma 2.4)}$$
$$= \frac{q + 2\sqrt{1-q^2}}{q} \max\{w_q\left(((T_1-T_2)(T_1+T_2))^n\right), w_q\left(((T_1+T_2)(T_1-T_2))^n\right)\}.$$

This proves the theorem. □

## 3. Bounds for $q$-numerical radius via Buzano inequality

In this section, we have obtained some bounds for $q$-numerical radius of operators using Buzano inequality. Further, we have investigated the $q$-numerical radius inequalities for commutators of positive operators.

**Lemma 3.1.** *Let $a, b \in \mathcal{L}(\mathcal{H})$ and $t \in \mathbb{R}$, then*

$$\|a\|^2\|b\|^2 - |\langle a, b\rangle|^2 \leq \|a\|^2\|b - ta\|^2 \tag{3.1}$$

*and*

$$\|a\|\|b\| + |\langle a, b\rangle| \geq 2|\langle a, c\rangle\langle c, b\rangle|. \tag{3.2}$$

**Theorem 3.2.** *Let $T \in \mathcal{L}(\mathcal{H})$. Then*

$$w_q^2(T) - w(T^2) \leq \frac{1}{|q|^2}\|\lambda T - T^*\|^2. \tag{3.3}$$

*Proof.* Let $b = T^*$, $a = \lambda Tx$ and $t = 1$, where $x \in \mathcal{H}, \|x\| = 1$ in (3.1), we have

$$|\lambda|^2\|Tx\|^2\|T^*x\|^2 \leq |\langle \lambda Tx, T^*x\rangle|^2 + \|\lambda Tx\|^2\|T^*x - \lambda Tx\|^2$$
$$= |\langle \lambda T^2x, x\rangle|^2 + |\lambda|^2\|Tx\|^2\|T^*x - \lambda Tx\|^2$$
$$= |\lambda|^2|\langle T^2x, x\rangle|^2 + |\lambda|^2\|Tx\|^2\|\lambda Tx - T^*x\|^2.$$

Taking supremum over $x \in \mathcal{H}$, $\|x\| = 1$, we get

$$\|T\|^2 \|T^*\|^2 \leq w^2(T^2) + \|T\|^2 \|\lambda T - T^*\|^2 \qquad (3.4)$$

on the other hand put $a = Ty$, $b = T^*y$, $e = x$, in (3.2), we have

$$\|Ty\| \|T^*y\| + |\langle Ty, T^*y \rangle| \geq 2|\langle Ty, x \rangle \langle x, T^*y \rangle|.$$

Equivalently, it can be written as

$$2|\langle Ty, x \rangle \langle x, T^*y \rangle| \leq |\langle T^2 y, y \rangle| + \|Ty\| \|T^*y\|.$$

So, we get

$$2|\langle T^*x, y \rangle \langle Tx, y \rangle| \leq |\langle T^2 y, y \rangle| + \|Ty\| \|T^*y\|,$$

which implies

$$w_q(T^*) w_q(T) \leq |\langle T^2 y, y \rangle| + \|Ty\| \|T^*y\|.$$

Since, $w_q(T) = w_q(T^*)$ and $0 < q \leq 1$, we have

$$2 w_q^2(T) \leq w(T^2) + \|T\| \|T^*\|$$

Using (3.4), we have

$$(2 w_q^2(T) - w(T^2))^2 \leq w^2(T^2) + \|T^2\| \|\lambda T - T^*\|^2,$$

which is same as

$$4 w_q^4(T) + w^2(T^2) - 4 w_q^2(T) w(T^2) \leq w^2(T^2) + \|T^2\| \|\lambda T - T^*\|^2.$$

This implies that

$$4 w_q^2(T)(w_q^2(T) - w(T^2)) \leq \frac{4}{|q|^2} w_q^2(T) \|\lambda T - T^*\|^2, \text{ by Lemma 1.13}$$

which implies

$$w_q^2(T) - w(T^2) \leq \frac{1}{|q|^2} \|\lambda T - T^*\|^2.$$

$\square$

**Remark 3.3.** *Taking $q \to 1$ in Theorem 3.2, we have*

$$w^2(T) - w(T^2) \leq \|\lambda T - T^*\|^2,$$

*which is Theorem 2.1 of [27].*

**Theorem 3.4.** *Let $T, S \in \mathcal{L}(\mathcal{H})$. Then*

$$w_q\left(\begin{bmatrix} O & T \\ S & O \end{bmatrix}\right) \leq \sqrt{\max(w(TS), w(ST)) + \frac{1}{|q|^2}\|T \pm S^*\|^2}.$$

*Proof.* Using (3.3), we have

$$w_q^2\left(\begin{bmatrix} O & T \\ S & O \end{bmatrix}\right) \leq w\left(\begin{bmatrix} O & T \\ S & O \end{bmatrix}^2\right) + \frac{1}{|q|^2}\left\|\begin{bmatrix} O & T \\ S & O \end{bmatrix} \pm \begin{bmatrix} O & S^* \\ T^* & O \end{bmatrix}\right\|^2$$

$$= w\left(\begin{bmatrix} TS & O \\ O & ST \end{bmatrix}\right) + \frac{1}{|q|^2}\left\|\begin{bmatrix} O & T \pm S^* \\ S \pm T^* & O \end{bmatrix}\right\|^2.$$

as $\|T \pm S^*\| = \|S \pm T^*\|$ and using equalities (1.10) and (1.11), we have

$$w_q^2\left(\begin{bmatrix} O & T \\ S & O \end{bmatrix}\right) \leq \max(w(TS), w(ST)) + \frac{1}{|q|^2}\|T \pm S^*\|^2.$$

□

**Remark 3.5.** *Taking $q \to 1$ in Theorem 3.4 we have*

$$w^2\left(\begin{bmatrix} O & T \\ S & O \end{bmatrix}\right) \leq \max(w(TS), w(ST)) + \|T \pm S^*\|^2,$$

*which is Theorem 2.2 of [27].*

**Lemma 3.6.** *[1, Theorem 5.5] If $T, S \in \mathcal{L}(\mathcal{H})$ and $q \in (0, 1]$, then*

$$w_q\left(\begin{bmatrix} O & T \\ S & O \end{bmatrix}\right) \leq \frac{q + 2\sqrt{1-q^2}}{2q}(w_q(T+S) + w_q(T-S)). \tag{3.5}$$

**Lemma 3.7.** *[1, Lemma 5.2] If $A, B \in \mathcal{L}(\mathcal{H})$, $q \in \overline{\mathbb{D}}$, then*

$$w_q\left(\begin{bmatrix} O & A \\ B & O \end{bmatrix}\right) = w_q\left(\begin{bmatrix} O & B \\ A & O \end{bmatrix}\right). \tag{3.6}$$

*and, by using (3.5),*

$$w_q(A) \leq w_q\left(\begin{bmatrix} O & A \\ A & O \end{bmatrix}\right) \leq \frac{q + 2\sqrt{1-q^2}}{q}w_q(A). \tag{3.7}$$

**Theorem 3.8.** *If* $P, Q, R \in \mathcal{L}(\mathcal{H})$. *Then*

$$w_q(PR \pm R^*Q) \le 2\sqrt{\max(w(PRR^*Q), w(R^*QPR)) + \frac{1}{|q|^2}\|PR \pm Q^*R\|^2} \qquad (3.8)$$

*Proof.* Using (3.6) and (3.7), we have

$$w_q\left(\begin{bmatrix} O & PR+R^*Q \\ PR+R^*Q & O \end{bmatrix}\right) = w_q(PR+R^*Q)$$

$$= w_q\left(\begin{bmatrix} O & R^*Q \\ PR & O \end{bmatrix} + \begin{bmatrix} O & PR \\ R^*Q & O \end{bmatrix}\right)$$

$$\le w_q\left(\begin{bmatrix} O & R^*Q \\ PR & O \end{bmatrix}\right) + w_q\left(\begin{bmatrix} O & PR \\ R^*Q & O \end{bmatrix}\right)$$

$$= w_q\left(\begin{bmatrix} O & PR \\ R^*Q & O \end{bmatrix}\right) + w_q\left(\begin{bmatrix} O & PR \\ R^*Q & O \end{bmatrix}\right)$$

$$= 2w_q\left(\begin{bmatrix} O & PR \\ R^*Q & O \end{bmatrix}\right).$$

So, we have

$$w_q(PR + R^*Q) \le 2w_q\left(\begin{bmatrix} O & PR \\ R^*Q & O \end{bmatrix}\right). \qquad (3.9)$$

By Theorem 3.4, we have

$$w_q(PR + R^*Q) \le 2\sqrt{\max(w(PRR^*Q), w(R^*QPR)) + \frac{1}{|q|^2}\|PR \pm Q^*R\|^2}.$$

Replacing $R$ by $iR$ in the above inequality, we have

$$w_q(PR - R^*Q) \le 2\sqrt{\max(w(PRR^*Q), w(R^*QPR)) + \frac{1}{|q|^2}\|PR \pm Q^*R\|^2}.$$

$\square$

**Corollary 3.9.** *Let* $P, Q, U \in \mathcal{L}(\mathcal{H})$ *and* $U$ *is a unitary operator, then*

$$w_q(PU \pm U^*Q) \le 2\sqrt{\max(w(PQ), w(QP)) + \frac{1}{|q|^2}\|P \pm Q^*\|^2}.$$

*In particular,*

$$w_q(PU \pm U^*Q) \le 2\sqrt{w(P^2) + \frac{1}{|q|^2}\|P \pm P^*\|^2}. \qquad (3.10)$$

*Proof.* Put $R = U$ in (3.8), we have

$$w_q(PU \pm U^*Q) \le 2\sqrt{\max(w(PUU^*Q), w(U^*QPU)) + \frac{1}{|q|^2}\|PU \pm Q^*U\|^2}. \quad (3.11)$$

Since $U$ is unitary and $w(U^*QPU) = w(QP)$, now the required result follows from (3.11). □

**Corollary 3.10.** *Let $P, Q, R, S \in \mathcal{L}(\mathcal{H})$ and $T = \begin{bmatrix} P & Q \\ R & S \end{bmatrix}$, then $\max\{w_q(P), w_q(S)\} \le \sqrt{w(T^2) \pm \frac{1}{|q|^2}\|T \pm T^*\|^2}$.*

*Proof.* By using Lemma 1.5

$$\max\{w_q(P), w_q(S)\} \le w_q\left(\begin{bmatrix} P & O \\ O & S \end{bmatrix}\right).$$ Since $S$ is arbitrary, we have $\max\{w_q(P), w_q(S)\} \le w_q\left(\begin{bmatrix} P & O \\ O & -S \end{bmatrix}\right)$. Let $U = \begin{bmatrix} I & O \\ O & -I \end{bmatrix}$ be a unitary operator on $\mathcal{H} \oplus \mathcal{H}$. Then,

$$\begin{aligned}
TU + U^*T &= \begin{bmatrix} P & Q \\ R & S \end{bmatrix}\begin{bmatrix} I & O \\ O & -I \end{bmatrix} + \begin{bmatrix} I & O \\ O & -I \end{bmatrix}\begin{bmatrix} P & Q \\ R & S \end{bmatrix} \\
&= \begin{bmatrix} P & -Q \\ R & -S \end{bmatrix} + \begin{bmatrix} P & Q \\ -R & -S \end{bmatrix} \\
&= \begin{bmatrix} 2P & O \\ O & -2S \end{bmatrix} \\
&= 2\begin{bmatrix} P & O \\ O & -S \end{bmatrix}.
\end{aligned}$$

Thus using (3.10), we have

$$\begin{aligned}
\max\{w_q(P), w_q(S)\} &\le w_q(TU + U^*T) \\
&\le \sqrt{w(T^2) \pm \frac{1}{|q|^2}\|T \pm T^*\|^2}.
\end{aligned}$$

□

**Theorem 3.11.** *Let $T, P \in \mathcal{L}(\mathcal{H})$ such that $P$ is a projection. Then*

$$w_q(TP - PT) \le \sqrt{w(T^2) + \frac{1}{|q|^2}\|T \pm T^*\|^2}. \quad (3.12)$$

*Proof.* Let us use the decomposition $\mathcal{H} = ran(P) \oplus ker(P)$ and using (1.9) we can represent $P$ as the form $P = \begin{bmatrix} I_1 & O \\ O & O \end{bmatrix}$, where $I_1$ is the identity operator on $ran(P)$. With respect to this decomposition, $T$ can be written as $T = \begin{bmatrix} T_{11} & T_{12} \\ T_{21} & T_{22} \end{bmatrix}$. Then

$$PT - TP = \begin{bmatrix} I_1 & O \\ O & O \end{bmatrix} \begin{bmatrix} T_{11} & T_{12} \\ T_{21} & T_{22} \end{bmatrix} - \begin{bmatrix} T_{11} & T_{12} \\ T_{21} & T_{22} \end{bmatrix} \begin{bmatrix} I_1 & O \\ O & O \end{bmatrix}$$

$$= \begin{bmatrix} T_{11} & T_{12} \\ O & O \end{bmatrix} - \begin{bmatrix} T_{11} & O \\ T_{21} & O \end{bmatrix}$$

$$= \begin{bmatrix} O & T_{12} \\ -T_{21} & O \end{bmatrix}.$$

Suppose $I_2$ is the identity operator on $ker(P)$ and if $U = \begin{bmatrix} I_1 & O \\ O & -I_2 \end{bmatrix}$, then $U$ is unitary and

$$U^*T - TU = \begin{bmatrix} I_1 & O \\ O & -I_2 \end{bmatrix} \begin{bmatrix} T_{11} & T_{12} \\ T_{21} & T_{22} \end{bmatrix} - \begin{bmatrix} T_{11} & T_{12} \\ T_{21} & T_{22} \end{bmatrix} \begin{bmatrix} I_1 & O \\ O & -I_2 \end{bmatrix}$$

$$= \begin{bmatrix} T_{11} & T_{12} \\ -T_{21} & -T_{22} \end{bmatrix} - \begin{bmatrix} T_{11} & -T_{12} \\ T_{21} & -T_{22} \end{bmatrix}$$

$$= 2 \begin{bmatrix} O & T_{12} \\ -T_{21} & O \end{bmatrix}.$$

So,
$$\frac{1}{2}(U^*T - TU) = \begin{bmatrix} O & T_{12} \\ -T_{21} & O \end{bmatrix}.$$

Similarly,
$$TU - U^*T = \begin{bmatrix} T_{11} & -T_{12} \\ T_{21} & -T_{22} \end{bmatrix} - \begin{bmatrix} T_{11} & T_{12} \\ -T_{21} & -T_{22} \end{bmatrix} = 2 \begin{bmatrix} O & -T_{12} \\ T_{21} & O \end{bmatrix}.$$

So, we have $\frac{1}{2}(TU - U^*T) = \begin{bmatrix} O & -T_{12} \\ T_{21} & O \end{bmatrix}$ and $TP - PT = \begin{bmatrix} T_{11} & O \\ T_{21} & O \end{bmatrix} - \begin{bmatrix} T_{11} & T_{12} \\ O & O \end{bmatrix}$.

This implies
$$TP - PT = \begin{bmatrix} O & -T_{12} \\ T_{21} & O \end{bmatrix}.$$

So,
$$w_q(TP - PT) = w_q\left(\begin{bmatrix} O & -T_{12} \\ T_{21} & O \end{bmatrix}\right) = \frac{1}{2}w_q(TU - U^*T).$$

Now, using (3.10), we have $w_q(TP - PT) \leq \sqrt{w(T^2) + \frac{1}{|q|^2}\|T \pm T^*\|^2}$. $\square$

**Theorem 3.12.** *Let $T, X \in \mathcal{L}(\mathcal{H})$ such that $T$ is positive. Then*
$$w_q(TX - XT) \leq \|T\|\sqrt{w(X^2) + \frac{1}{|q|^2}\|X \pm X^*\|^2}.$$

*Proof.* Since $T\sqrt{T - T^2} = (\sqrt{T - T^2})T$, so $P = \begin{bmatrix} T & \sqrt{T - T^2} \\ \sqrt{T - T^2} & I - T \end{bmatrix}$ is a projection on $\mathcal{H} \oplus \mathcal{H}$.

Let $S = \begin{bmatrix} X & O \\ O & O \end{bmatrix}$, then

$$PS - SP = \begin{bmatrix} T & \sqrt{T - T^2} \\ \sqrt{T - T^2} & I - T \end{bmatrix}\begin{bmatrix} X & O \\ O & O \end{bmatrix} - \begin{bmatrix} X & O \\ O & O \end{bmatrix}\begin{bmatrix} T & \sqrt{T - T^2} \\ \sqrt{T - T^2} & I - T \end{bmatrix}$$
$$= \begin{bmatrix} TX & O \\ \sqrt{T - T^2}X & O \end{bmatrix} - \begin{bmatrix} XT & X\sqrt{T - T^2} \\ O & O \end{bmatrix}$$
$$= \begin{bmatrix} TX - XT & -X\sqrt{T - T^2} \\ (\sqrt{T - T^2})X & O \end{bmatrix}.$$

Now, using (3.12),
$$w_q(SP - PS) \leq \sqrt{w(S^2) + \frac{1}{|q|^2}\|S \pm S^*\|^2}.$$

Let $Q = \begin{bmatrix} I & O \\ O & O \end{bmatrix}$, then

$$Q(PS - SP)Q^* = \begin{bmatrix} I & O \\ O & O \end{bmatrix}\begin{bmatrix} TX - XT & -X\sqrt{T - T^2} \\ \sqrt{T - T^2}X & O \end{bmatrix}\begin{bmatrix} I & O \\ O & O \end{bmatrix}$$
$$= \begin{bmatrix} I & O \\ O & O \end{bmatrix}\begin{bmatrix} TX - XT & O \\ \sqrt{T - T^2}X & O \end{bmatrix}$$
$$= \begin{bmatrix} TX - XT & O \\ OX & O \end{bmatrix}.$$

Using Lemma 1.5 we have,

$$w_q(TX - XT) \leq w_q\left(\begin{bmatrix} TX - XT & O \\ O & O \end{bmatrix}\right) = w_q(Q(PS - SP)Q^*).$$

Since $w_q((Q(PS - SP)Q^*) \leq w_q(PS - SP)$. So,

$$w_q(TX - XT) \leq w_q\left(\begin{bmatrix} TX - XT & O \\ O & O \end{bmatrix}\right) = w_q(Q(PS - SP)Q^*) \leq w_q(PS - SP).$$

This implies,

$$w_q(TX - XT) \leq w_q(PS - SP)$$
$$\leq \sqrt{w(S^2) + \frac{1}{|q|^2}\|S \pm S^*\|^2} \quad \text{by (3.12)}$$
$$= \sqrt{w(X^2) + \frac{1}{|q|^2}\|X \pm X^*\|^2} \quad \text{by (1.10) and (1.11)}.$$

Now, since $T$ is a positive operator, $w_q\left(\frac{T}{\|T\|}X - X\frac{T}{\|T\|}\right) \leq \sqrt{w(X^2) + \frac{1}{|q|^2}\|X \pm X^*\|^2}$. Hence,

$$w_q(TX - XT) \leq \|T\|\sqrt{w(X^2) + \frac{1}{|q|^2}\|X \pm X^*\|^2}.$$

□

**Data availability:** The author declare that data sharing is not applicable to this article as no datasets were generated or analysed during the current study.

**Conflict of interest:** The author declare that there is no conflict of interest.